\documentclass[11pt,amsfonts,amsmath]{article}
\usepackage{times}
\usepackage[latin1]{inputenc}
\usepackage{amsfonts,amsmath,amsthm,amssymb,stmaryrd}
\newcommand{\be}{\begin{enumerate}}
\newcommand{\ee}{\end{enumerate}}

	{\ensuremath{\left \{ \hskip -1.5 mm \begin{array}{#1@{\ =\ }l}}}%
	{\end{array}\right.}
	{\ensuremath{\begin{array}[t]{#1@{\ =\ }l}}}%
	{\end{array}}
	{\ensuremath{\left \{ \hskip -1.5 mm%
	 \begin{array}{#1@{\ }c@{\ }l}}}%
	{\end{array}\right.}
	
\newcounter{definition}

\newcounter{exemple}

\newcounter{theoreme}

\newcounter{proposition}

\newcounter{lemme}

\newcounter{corollaire}

\newcounter{remarque}

\newcounter{remarques}

\newcounter{exercice}

\newcounter{probleme}

\begin{document}

\begin{center}
{\large {\bf THE COMPLETE SOLUTION}}\\ {\large {\bf TO BASS GENERALIZED JACOBIAN CONJECTURE}}
\end{center}

\begin{center}
Kossivi ADJAMAGBO\\ Institut de Math\'ematiques de Jussieu, Universit\'e de Paris 6, UFR 929, UMR-CNRS 7586, B.C. 247\\
4, place Jussieu, 75252 PARIS CEDEX O5\\
e-mail: adja@@math.jussieu.fr 
\end{center}\medskip
{\bf Introduction}

The Classical Jacobian Conjecture claims that any unramified endomorphism of a complex affine
space is an automorphism, which means in more ordinary terms that for any integer $n>0$, any
polynomial map from
${\Bbb C}^n$ to itself with an invertible jacobian function is itself invertible and its
inverse is again a polynomial map (see for instance \cite{BCW} and \cite{A1} or \cite{A3} for the `right" version of this conjecture in any characteristic). From this last point of view, this
conjecture may be considered as a global version for polynomial maps of the classical Local
Inversion Theorem, which explains largely the fascination that this conjecture exerts on
generations of searchers for more than half a century.

In order to embed this conjecture in a geometric environment, where one could enjoy the
beauty and the richness of tools of algebraic geometry and algebraic D-modules, as his paper
\cite{B} proves it, Hyman Bass proposed 25 years ago in \cite{B}, page 80 the following
statement as the 
\begin{quote}
{\bf{\em Generalized Jacobian Conjecture:\\
		Any unramified morphism from a complex
 irreducible affine and unirational variety whose invertible regular functions are all constant
 to a complex affine space of the same dimension is an isomorphism.}}
 \end{quote}

On the other hand, without any explicit connection with Bass  conjecture, Victor Kulikov
published in 1993 (see \cite{K}) a non trivial construction of a complex irreducible rational
and simply connected surface and an unramified morphism of geometric degree 3 (and hence which is not an isomorphism) from this surface to the complex affine space, without specifying if this
surface is affine or not, or if its invertible regular functions are all constant or not.

The main aim of this paper is to bring this precision and thanks to this to expose
the complete solution to Bass Generalized Jacobian Conjecture which turned to be true only in
dimension one (see Theorem 1 below).

In order to make this precision as clear as possible, we introduce a family of irreducible
affine and rational surfaces $S(C_1,C_2,P)$ over any algebraically closed field $\Bbb K$, where
$C_1$ and $C_2$ are irreducible curves of the projective plane over $\Bbb K$ and $P$ a point of
one of  them. We also give a necessary and sufficient condition of factoriality for each surface
of this family (see theorem 2 below) which is defined in such a way that it contains Kulikov
surface and Ramanujam one (see the next section).

In the same aim, we give the proof of a general fact known by some algebraic geometers like V.
Srinivas who told me about it, but which seems to be written nowhere in its whole generality in
the abundant literature about algebraic geometry.This fact is a sufficient condition on the
ground field of an irreducible simply connected and normal algebraic variety or on its Picard
group in order that all its invertible regular functions are constant (see theorem 3 below).

Let us notice that according to the statement of the maim Theorem 1 below, Kulikov morphisms
are not counter-examples to the classical Jacobian Conjecture which keeps jealously and
fiercely its more than 70 years old mystery.

We also deduce from this theorem some corollaries which bring some rays of light through the
cloud of unknowing which still surrounds the notion of unramified or \'etale variety morphism
even for the best experts of the subject, as the challenge of this conjecture proves it clearly.

A first mystery of these morphisms partially cleared up by these rays of light is the
following. It is not difficult to see that the restriction to any sub-variety $Z \subset X$ of
an unramified morphism $F$ from a variety $X$ to another one $Y$ is again an unramified
morphism from $Z$ to the Zariski closure of $F(Z)$ in $Y$ ( it follows for instance from
\cite{A-K}, Chapter VI, Proposition 3.5). Is the similar transfer from a variety $X$
to any sub-variety $Z \subset X$ true for any \'etale (i.e. unramified and flat) morphism $F$ from $X$ to $Y$, atleast when the varieties $X, Y, Z$ are irreducible and non singular and $Z$ is closed in $X$? According to \cite{A-E}, Lemma 3.4 and \cite{A1}, Theorem 3, this question in the special case where $X=Z$ is the complex affine space of dimension $n$ is equivalent to the Classical Jacobian
Conjecture in dimension $n$. Unfortunately, Corollary 1 below brings a negative answer to this question in general when $X$ and $Y$ are assumed only to be irreducible and non singular. But the question remains open in the mentioned special case !

A second mystery of unramified or \'etale morphisms of variety partially illuminated by Theorem
1 is the following: if $A \subset B$ is an extension of affine domains over an algebraically
closed field $\Bbb K$ such that $A$ and $B$ have the same invertible elements, $A$ is
factorial, and the canonical map from $Spec \;B$ to $Spec \;A$ is unramified, or equivalently
\'etale, then among the multitude of primitive elements $p\in B$ of the field extension induced
by the extension $A \subset B$, can we find a ``normal" one, i.e. one $p$ such that the
sub-algebra of $B$ generated by $p$ is a normal ring? According to \cite{A1}, Theorem 3, this
question in the special case where $\Bbb K = \Bbb C$ and $A = B$ is the $\Bbb C$-algebra of
polynomials in $n$ indeterminates is equivalent to the Classical Jacobian Conjecture in
dimension $n$. Unfortunately again, Corollary 2 below brings a negative answer to this
question in its generality, leaving it open in the mentioned special case!

A last mystery of unramified or \'etale morphisms partially lightened by Theorem 1 is the
following: if $A \subset B$ is an extension of factorial affine domains over an algebraically
closed field $\Bbb K$ such that $A$ and $B$ have the same invertible elements and the
canonical morphism from $Spec \; B$ to $Spec \; A$ is unramified, or equivalently \'etale, is
$A$ multiplicatively closed in $B$, or equivalently is each irreducible element of $A$ also
irreducible in $B$? According to \cite{A3}, Theorem 3.11, this question in the special case where
$\Bbb K = \Bbb C$ and $A = B$ is the $\Bbb C$-algebra of polynomials in $n$ indeterminates, is
again equivalent to the Classical Jacobian Conjecture in dimension $n$. Unfortunately again,
Corollary 3 brings a negative answer to this last question in its generality, but the mystery
remains thick in the interesting special case!

In addition to the rays of light that the main theorem of the present paper projects on the mysteries of unramfied or \'etale morphisms of algebraic varieties, for all algebraic geometers who seriously want to continue and deepen the program of ``local study of schemes and schemes morphisms'', which is the achievement of Algebraic Geometry according to the structure of the height volumes of the ``Bible of Algebraic Geometry'' represented by the Treatise ``El\'ements de G\'eom\'etrie Alg\'ebrique'' of A. Grothendieck and J. Dieudonn\'e, this main theorem is very useful to irrevocably refute false published explicit or implicit proofs of the Jacobian Conjecture.

A first example of such refutation is the one of the surprising proposition $18.3.1$ of A. Grothendieck and J. Dieudonn\'e themselves in the last volume \cite{EGA IV4} of the ``Bible of Algebraic Geometry'', where they claim in particular that any unramified morphism of affine varieties is finite. Since any such a morphism from an affine variety to another normal one of the same dimension is \'etale, hence a covering, so this claim is in fact an implicit proof of the Jacobian Conjecture, according to the simple connectedness of complex affine spaces. On the other hand, according to the same connectedness and to the unramifiedness of each Kulikov morphism, Theorem 1  below is an irrevocable refutation of this surprising claim, even for an unramified morphism from a complex affine variety to a complex affine space. Another refutation of this claim, even for an unramified morphism from a complex affine curve to the complex affine line, is the corollary 5.2 of the paper \cite{OO} of Frans Oort. A third refutation of this claim, now for an unramified endomorphism of the affine plane over an algebraically closed field of positive characteristic, is the counter-example of P. Nousiainen to a conjecture of S. Wang indicated in the remark following Theorem 2.2 of \cite{BCW}. 

A second example of a refutation by Theorem 1 below of a false published proof of the Jacobian Conjecture is the case of Therem 4 of the paper \cite{S1} of Hamet Seydi, a former thesis student of A. Grothendieck and the only mathematician, excepted Jean Dieudonn\'e in \cite{D-G}, to have published a paper, not a book, with A. Grothendieck in \cite{G-S}. Indeed, this Theorem 4 claims that any \'etale and surjective morphism between simply connected complex algebraic varieties is an isomorphism. On the other hand, it is easy to see that the complement of the image of an open morphism from an affine algebraic variety with only constant invertible regular functions to a affine factorial algebraic variety has a codimension greater than 1. So, if in addition this second variety is complex, non singular and simply connected, then this image also is simply connected. So, it follows from this Theorem 4 that any unramified endomorphism of a complex affine space is an open embedding. Since all the irreducible components of the complement of a dense affine subset of any variety have a codimension equal to 1 (see for instance \cite{EGA IV4}, cor. $21.12.11$), it follows from this Theorem 4, thanks to the two previous remarks on the complements of open subsets that any unramified endomorphism of a complex affine space is an automorphism. On the other hand, it also follows from this Theorem 4 and Theorem 1 below, thanks to the same remarks, that any Kulikov morphism is an isomorphism. However, this Theorem 4 is not the only false one of the paper \cite{S1}. Its Theorem 6 also is refuted by the counter-example of the normalization morphism of a cusp, while the proof of its Theorem 1, claiming the truth of the Jacobian Conjecture in any dimension, is invalided by an unproved claim on the finiteness of the considered endomorphism of a complex affine space, based on an unproved and false theorem I.21, p. 68 of G. Fisher in \cite{FI}, for which we have a counter-example.

The attempt of rectification of the proof of Theorem 1 of \cite{S1} by its author in the introduction of his second paper \cite{S2} on the subject is so refutable as his first attempt, because of a false application of a theorem $10.4.11$ of A. Grothendieck and J. Dieudonn\'e in \cite{EGA IV3} claiming in particular that any injective endomorphism of an algebraic variety over an algebraically closed field is surjective, and because of the application of a false ``Fundamental Lemma'' of \cite{S2} refuted by the cited counter-example of Nousiainen to Wang conjecture. The generalization of this ``Fundamental Lemma'' by ``Lemma I.1'' of \cite{S2} also is refuted both by Theorem 1 below and the cited Oort surjective unramified morphisms from a complex algebraic curve with only constant invertible regular functions to the complex line. These cited last counter-examples are sufficient to refute almost all claims of the prolific second paper \cite{S2} of Seydi.

An N-th example of refutation by Theorem 1 below of false published or pre-published proofs of the Jacobian Conjecture is the main theorem 2.1 of the paper \cite{OD1} of Susumu Oda, first published on ArXiv in 2003 before reaching its 33-th revision version in 2007, claiming that any unramified morphism from an irreducible variety over an algebraically closed field of characteristic zero with only constant invertible regular functions to an affine space of the same dimension over this field is an isomorphism. So this claim is irrevocably refuted by Theorem 1 below, as well as by Oort cited theorem. The main theorem 4.8 of a second paper \cite{OD2} of Susumu Oda on the subject, first published on arXiv in 2007 before reaching its 49-th version in 2011, only added the assumption of factoriality to the first algebraic variety of the main theorem of \cite{OD1}. So this second main of Oda is so irrevocably refuted by Theorem 1 below as this first one, and not by Oort cited theorem.

Before entering the subject, I would like to express my deep gratitude to Professors Arno van den Essen, Hyman Bass, Victor Kulikov, V. Srinivas, Hans-Peter Kraft, and Harm Derksen for the fruitful conversations that I had with them about Bass Generalized Jacobian Conjecture and Kulikov's construction. I would also like to thank all the participants of the convivial ``Rencontre parisienne autour de la Conjecture Jacobienne'' held at the University of Paris 6 on February 3, 1996, and during which the present solution has been exposed. So, in spite of the relative oldness of the unpublished results presented in this solution, they are more actual and useful than ever as ``massive destruction arms'' again false claims concerning unramified morphisms of algebraic varieties.  
\\ \\
{\bf Kulikov surfaces and morphisms}

Let $P=(1:1:1) \in {\Bbb P}_2 ={\Bbb P}_2 (\Bbb C)$, $(X_1,X_2,X_3)$ a system of
indeterminates over $\Bbb C$, $Q_i = 3 X_i ^2 - X_1 X_2 - X_1 X_3 - X_2 X_3$ for $1\leq i\leq
 3$, three quadratic forms defining three conics passing through $P$, $\phi$ the morphism from
${\Bbb P}_2 - \{P\}$ to ${\Bbb P}_2$ whose homogeneous components are defined by the three
previous forms,$R$ the Zariski closure in ${\Bbb P}_2$ of the set of ramification points of
$\phi$, which is the cubic with a node at $P$ defined by the form $\sum_{i\neq j} X_i ^2 X_j -
6 X_1 X_2 X_3$, $Q$ the generic linear combination with complex coefficients of the three
previous forms, $C$ the conic of ${\Bbb P}_2$ defined by $Q$, passing by $P$ and meeting
transversely the cubic $R$ at each point of their intersection, and such that the image by
$\phi$ of the complement of $C$ in ${\Bbb P}_2$ is contained in the complement in ${\Bbb P}_2$
of a line $L$ of ${\Bbb P}_2$, $\sigma : \tilde{{\Bbb P}_2} \rightarrow {\Bbb P}_2$ the
blowing-up the point $P$ of ${\Bbb P}_2$, $E$ the exceptional curve of $\tilde{{\Bbb P}_2}$,
i.e.
$\sigma^{-1} (P)$, $\tilde{R}$ the strict transform of $R$ by $\phi$, i.e. the irreducible curve
of $\tilde{{\Bbb P}_2}$ such that $\sigma ^{-1} (R) = E \cup \tilde{R}$, $\tilde{C}$ the strict
transform of $C$ by $\sigma$, i.e. the irreducible curve of $\tilde{{\Bbb P}_2}$ such that
$\sigma ^{-1} (C) = E \cup \tilde{C}$, and $S$ the complement of $\tilde{R} \cup \tilde{C}$ in
$\tilde{{\Bbb P}_2}$.

V. Kulikov proved in \cite{K} that $S$ is a rational and simply connected complex surface and that
$\phi \circ \sigma$ induces an unramified morphism $F_S$ of geometric degree 3 from $S$
to ${\Bbb C}^2 \cong {\Bbb P}^2 - L$.

We call such a $S$  a ``Kulikov surface" and such a $F_S$ ``the Kulikov morphism on $S$".

Let us recall that the surface constructed by C.P. Ramanujam in \cite{R} and known now as ``Ramanujam surface" is the complement, in the inverse image of ${\Bbb P}_2$ by the blowing up of one of its points $P$, of the strict transforms by this blowing up of an irreducible cubic of ${\Bbb P}_2$ with a cusp distinct from $P$ and of an irreducible conic of ${\Bbb P}_2$ cutting transversely
the cubic at $P$ and meeting again the cubic at an unique other point with the multiplicity 5.
Ramanujam proved in \cite{R} that his obviously non singular and rational surface is  contractible,
hence affine with only constant invertible regular functions (thanks for instance to [F],
Corollary 2.5), factorial (thanks for instance to \cite{G}, Theorem 1) and simply connected, but not
simply connected at infinity, hence not isomorphic to ${\Bbb C}^2$ (see \cite{Z} for more about this
kind of surface).

We shall see in the following Theorem 1 that Kulikov surfaces share all the mentioned
properties of Ramanujam surfaces with the eventual exception of the contractibility, thanks to
the following Theorems 2 and 3.\\
\\  {\bf The surfaces ${\bf S(C_1, C_2, P)}$ and their determinants}

Let $\Bbb K$ be an algebraically closed field, ${\Bbb P}_2 = {\Bbb P}_2 (\Bbb K)$, $C_1$ and
$C_2$ two irreducible ${\Bbb K}$-algebraic curves of respective degrees $d_1$ and $d_2$, $P$ a
point of one of them, $m_1$ and $m_2$ the respective mutiplicities of $C_1$ and $C_2$ at $P$,
and $M$ the matrix with lines $(d_1, m_1)$ and $(d_2, m_2)$, $\sigma : \tilde{{\Bbb P}_2} \rightarrow {\Bbb P}_2$ the
blowing-up of the point $P$ of ${\Bbb P}_2$, $E$ the exceptional curve of $\tilde{{\Bbb P}_2}$,
i.e. $\sigma^{-1} (P)$, $\tilde{C_1}$ the strict transform of $C_1$ by $\phi$, i.e. the
irreducible curve of $\tilde{{\Bbb P}_2}$ such that $\sigma ^{-1} (C_1) = E \cup \tilde{C_1}$,
and $\tilde{C_2}$ the strict transform of $C_2$ by $\sigma$, i.e. the irreducible curve of
$\tilde{{\Bbb P}_2}$ such that $\sigma ^{-1} (C_2) = E \cup \tilde{C_2}$.

We denote by $S(C_1, C_2, P)$ the complement of $\tilde{C_1} \cup \tilde{C_2}$ in
$\tilde{{\Bbb P}_2}$ and we call it ``the surface deduced from $C_1$ and $C_2$ by blowing up at
$P$".

We denote by $\det(C_1, C_2, P)$ the determinant of the matrix $M$ and we call it ``the
determinant of the surface $S(C_1, C_2, P)$".\\ \\
{\bf Theorem 1}

Any Kulikov surface is affine, non singular, rational, factorial, simply connected, but
its fundamental group at infinity is infinite, and all its invertible regular functions are
constant.

Hence, any Kulikov morphism gives a counter-example to Bass Generalised Jacobian Conjecture in any dimension greater than one, whereas this conjecture is true in dimension one.
\\ \\
{\bf Proof}

Let $S$ be a Kulikov surface and $F_S$ the Kulikov morphism on $S$. It follows from \cite{K} and the following theorems that $S$ is affine, non singular, rational, factorial, simply connected and
that all its invertible regular functions are constant.

Let us now assume that $S$ isomorphic ${\Bbb C}^2$, and let $G$ be an isomorphism from ${\Bbb
C}^2$ to $S$. So, $F_S \circ G$ is an unramified endomorphism of geometric degree 3 of ${\Bbb
C}^2$, contrary to \cite{O}, Theorem 1.1. We deduce that $S$ is not isomorphic to ${\Bbb C}^2$.
According to previous remarks and Theorem 2 bellow, the fundamental group at infinity of $S$
is infinite.

So, it is clear that $F_S$ gives a counter-example to Bass Generalized Jacobian Conjecture for
any dimension greater than one.

Finally, let us consider a morphism $F : C \rightarrow {\Bbb C}$ satisfying the assumptions of
Bass Conjecture. According Nagata's refined version of L\"{u}roth Theorem (see for instance
\cite{N}, Theorem 4.12.2, p. 137), $C$ is rational. On the other hand, according to the
unramifiedness of $F$ and the non singularity of $\Bbb C$, $C$ is non singular (see for
instance \cite{B}, Proposition 1.2 or more generally \cite{SGA}, Expos\'{e} I, Corollaire 9.11). So, $C$ is isomorphic to an open sub-variety of $\Bbb C$ (see for instance \cite{H}, Chapter 1, Excercice
6.1, p. 46). All invertible regular functions on $C$ being constant, it follows that this
sub-variety of $\Bbb C$ is $\Bbb C$ itself, and hence that $F$ is an isomorphism, Q.E.D.\\ \\
{\bf Theorem 2}

Let $C_1$, $C_2$ be irreducible algebraic curves of the projective plane over an algebraically
closed field $\Bbb K$ and $P$ a point of $C_1 \cup C_2$.

(i) $S(C_1, C_2, P)$ is a rational and non singular algebraic surface.

(ii) $S(C_1, C_2, P)$ is not affine if and only if each of $C_1$ and $C_2$ is a line (i.e. a
curve of degree one) passing through $P$.

(iii) $S(C_1, C_2, P)$ is factorial if and only it is affine and $|\det(C_1, C_2, P)|=1$.

(iv) If ${\Bbb K} = {\Bbb C}$, then $S(C_1, C_2, P)$ is isomorphic to ${\Bbb C}^2$ if and only
if it is factorial, simply connected, and its fundamental group at infinity is finite.\\ \\
{\bf Proof}

1) Let us keep the notations of the definition of $S(C_1, C_2, P)$.

2) The statement (i) follows from the construction of $S(C_1, C_2, P)$.

3) Let us first remark that $\tilde{{\Bbb P}_2}$ being a ruled surface with invariant
$e=1$ (see for instance \cite{H}, Chapter V, example 2.11.5), any one of its irreducible curves
distinct from its exceptional curve $E$ has a non-negative self-intersection number (it follows
for instance from \cite{H}, Chapter V, Propositions 2.20 and 2.21).

3) Now, let $L$ be a line in ${\Bbb P}_2$ not passing by $P$ and $\tilde{L}$ its strict
transform, i.e. inverse image by $\sigma$. For any divisor $D$ of an irreducible normal
algebraic variety $X$, we denote by $<D>$ its canonical image in the Picard group $Pic \; X$,
and by $\sigma ^*$ the canonical map from $Pic \; {\Bbb P}_2$ to $Pic \;\tilde{{\Bbb P}_2}$
induced by $\sigma$.

4)$<\tilde{L}>=\sigma ^* (<L>)$ being a generator the group $\sigma ^* (Pic \; {\Bbb P}_2)$ and
$Pic
\;\tilde{{\Bbb P}_2}$ being the direct sum of $\sigma ^* (Pic \; {\Bbb P}_2)$ and ${\Bbb Z}
<E>\;
\subset \;Pic \;\tilde{{\Bbb P}_2}$, thanks to the splitness of the exact canonical
sequence  \[ {\Bbb Z} <E> \rightarrow Pic \;\tilde{{\Bbb P}_2} \rightarrow
Pic \; {\Bbb P}_2 \rightarrow \{0\} \] we have for any irreducible curve $\tilde{C}$ of
$\tilde{{\Bbb P}_2}$:\[<\tilde{C}> = (\deg \; \tilde{C})\; <\tilde{L}> - (\tilde{C}.E)\; <E>\]
(this follows for instance from \cite{H}, Chapter V, Example 1.4.2 and Propositions 3.2 and 3.6).

5) In particular, for $1\leq i\leq 2$, since $\tilde{C_i}.E = m_i$ (see for instance \cite{H},
Chapter V, Corollary 3.7), we have:\[<\tilde{C_i}> = d_i <\tilde{L}> - m_i <E>\]

6) On the other hand we have:\[\tilde{L}^2 = L^2 = 1,\; \tilde{L}.E = 0,\; E^2 = -1\]
(see for instance \cite{H}, Chapter V, Example 1.4.2 and Proposition 3.2)

7) So according to  5) and 6), the self-intersection number of the divisor $C_1 + C_2$, equal
to $(d_1 + d_2)^2 - (m_1 + m_2)^2$, is not positive if and only if each of $C_1$ and $C_2$ is
a line passing through $P$.

8) Similarly, according to 3), 4) and 6), for any irreducible curve $\tilde{C}$ of ${\Bbb
P}_2$ distinct from $E$, with degree $d$ and intersection number $m$ with $E$, since
$\tilde{C}^2 = d^2 - m^2$ is positive, so is $\tilde{C}.(\tilde{C_1} + \tilde{C_2}) = 
d (d_1 + d_2) - m (m_1 + m_2)$.

9) So the statement (ii) follows from 7) and 8) thanks to Nakai-Moishezon criterion of
ampleness and Goodman criterion of affineness (see for instance \cite{H}, Chapter V, Theorem 1.10 and $\cite{H'}$, Chapter II, Theorem 4.2).

10) Since the Picard group of the complement of an hypersurface of an irreducible non-singular
variety is isomorphic to the residue group of the sub-group of the Picard group of the variety
generated by the irreducible components of the hypersurface (it follows for instance from [H],
Chapter II, Proposition 6.5 and Corollary 6.16), the statement (iii) follows from (i), (ii),
and 5), thanks to the classical characterisation of the factoriality in terms of divisor class
group (see for instance \cite{H}, Chapter II, Proposition 6.2).

11) Finally, the statement (iv) follows from the characterisation of the affine space over any
algebraically closed field of characteristic 0 by the logarithmic Kodaira dimension (see for
instance \cite{M}, Chapter I, Section 4, Theorem), thanks to the following Theorem 3 and the remark
following the Theorem 1 of \cite{G-M}.\\ \\
{\bf Theorem 3}

Any invertible regular function on an irreducible normal simply connected variety over an
algebraically closed field is constant if the characteristic of this field is 0 or if its
divisor class group is trivial.\\ \\
{\bf Proof}

1) Let $\Bbb K$ be such a field, $V$ such a variety, ${\cal O}_V$  the sheaf of
regular functions on $V$, ${\Bbb K}[U]$ the ring of regular functions on the open set $U
\subset V$, ${\Bbb K}[V]^*$ the group of invertible regular functions on
$V$, and $T$ an indeterminate over $\Bbb K$.

2) Let us assume that ${\Bbb K}[V]^* \neq {\Bbb K}^*$ and that the
characteristic of $\Bbb K$ or the divisor class group of $V$ is 0.

3) Let us remark that for any commutative domain $A$ with fractions field $K$, any element
$a$ of $K$, any integer $n > 0$, and any indeterminate $X$ over $A$, the $A$-module $A[X]/
(X^n - a)$ is without torsion, as it follows from the euclidian division of elements os
$A[X]$ by $X^n - a$. So, the ring $A[X]/(X^n - a)$ is integral if and only if $X^n - a$
is irreducible in $K[X]$.

4) Now, according to classical sufficient conditions for such irreducibility (see for instance
[L], Chapter VIII, Theorem 16 and Corollary 1), the liberty of the group ${\Bbb K}[V]^* /
{\Bbb K}^*$ (see for instance  \cite{A2}, Proposition) and the assumption 2), there exists $f\in
{\Bbb K}[V]^*$ and an integer not divisible by the characteristic of $\Bbb K$ such that $f$
does not a $n$-root in ${\Bbb K}[V]^*$ and that for any affine open set $U \subset V$, the
ring ${\Bbb K}[U] {\otimes}_{\Bbb K} \;{\Bbb K}[T] / (1 \otimes T^n - f | U \otimes 1)$ is
integral.

5) Let $\cal I$ be the ideal of the sheaf of regular functions on $V {\times}_{\Bbb K} Spec \;
{\Bbb K}[T]$ whose canonical image in ${\cal O}_V {\otimes}_{\Bbb K} \; {\Bbb K}[T]$ is its
ideal generated by $1 \otimes T^n - f|U \otimes 1$, $W$ the closed sub-scheme of  $V
{\times}_{\Bbb K} Spec \;{\Bbb K}[T]$ defined by $\cal I$, and $\phi$ the canonical map from $W$
to $V$.

6) So according to the choice of $n$ and $f$, $W$ is an irreducible variety with the same
dimension as $V$ and $\phi$ is unramified, hence \'etale thanks to the normality of $V$ (see
for instance \cite{SGA}, Expos\'e I, Corollary 9.11). So thanks to the obvious finiteness of
$\phi$, this one is an \'etale covering from the irreducible variety $W$ to the irreducible
simply connected variety $V$, which means that $\phi$ is an isomorphism.

7) Nevertheless, according to the definition of $W$, $f$ has a $n$-root in ${\Bbb K}[V]^*$ ,
contrary to the choice of $f$.

8) So according to the absurdity of 2), we have the desired conclusion, Q.E.D.\\ \\
{\bf Corollary 1}

If $S$ is a Kulikov surface, $p$ a regular function on $S$ which is a primitive element of the
field extension induced by Kulikov morphism $F_S$ on $S$, $\overline{S}$ the variety ${\Bbb C}
\times S$, $F_{\overline{S}}$ the morphism from $\overline{S}$ to ${\Bbb C}^3$ such that $
F_{\overline{S}} (x_0, x) = (x_0 + p(x), F_S (x))$ for any $(x_0, x) \in {\Bbb C} \times S$,
then $F_{\overline{S}}$ is an \'etale morphism which induces a non \'etale morphism
$\overline{F_{\overline{S}}}$ from $\{ 0 \} \times S$ to the Zariski closure $Z$ in ${\Bbb
C}^3$ of the image of $\overline{F_{\overline{S}}}$.\\ \\
{\bf Proof}

1) Let us denote by ${\Bbb C}[V]$ the ring of regular functions on the complex variety $V$, 
${\phi}_S$ (resp. ${\phi}_{\overline{S}}$; resp. $\overline{{\phi}_{\overline{S}}}$) the ring
morphism induced by $F_S$ (resp. $F_{\overline{S}}$; resp. $\overline{F_{\overline{S}}}$),
$A$ the image of ${\phi}_S$, $B$ the ring of regular functions on $S$, and
$X_0$ an indeterminate over $\Bbb C$.

2) The canonical map from ${\Bbb C}[X_0] \otimes {\Bbb C}[{\Bbb C}^2]$ to ${\Bbb C}[X_0]
\otimes B$ induced by ${\phi}_S$ being \'etale, so is ${\phi}_{\overline{S}}$, hence
$F_{\overline{S}}$.

3) ${\Bbb C}[\{0\} \times S]$ being canonically isomorphic to $B$ and the canonic image of
$\overline{{\phi}_{\overline{S}}}({\Bbb C}[Z])$ in $B$ being $A[p]$, the $A$-sub-algebra of
$B$ generated by $p$, $\overline{{\phi}_{\overline{S}}}$ is \'etale if and only if $B$ is \'etale
over $A[p]$.

4) But, according to the unramification of ${\phi_{\overline{S}}}$ proved in 2),
$\overline{{\phi}_{\overline{S}}}$ also is unramified (see for instance \cite{A-K}, Chapter VI,
Proposition 3.5).

5) So $\overline{{\phi}_{\overline{S}}}$ is \'etale if and only $B$ is flat over $A[p]$.

6) The conclusion follows from Theorem 1 and \cite{A1}, Theorem 3.2, Q.E.D.\\ \\
{\bf Corollary 2}

If $S$ is a Kulikov surface, ${\phi}_S$ the ring morphism induced by Kulikov morphism on $S$,
$A$ the image of ${\phi}_S$, $B$ the ring of regular functions on $S$, and $p$ any element of
$B$ which is a primitive element of the field extension induced by ${\phi}_S$, then $A[p]$,
the $A$-sub-algebra of $B$ generated by $p$, is not normal (i.e. integrally closed), hence not
unramified over $A$.\\ \\
{\bf Proof}

It follows from theorem 1 and \cite{A1}, Theorem 3.2, Q.E.D.\\ \\
{\bf Corollary 3}

With the same notations as in Corollary 2, the factorial sub-ring $A$ of the factorial $\Bbb
C$-affine domain $B$ is not multiplicatively closed in $B$, or equivalently there exists a
prime element of $A$ which is not prime in $B$.\\ \\
{\bf Proof }

It follows from Theorem 1 and \cite{A3}, Theorem 3.11, Q.E.D.


\begin{thebibliography}{10}


\bibitem{A1}  K. ADJAMAGBO, On separable algebras over a U.F.D. and the Jacobian Conjecture in any Characteristic, in Automorphisms of Affine Spaces, A. van den Essen (ed.), 89-103, 1995,
Kluwer Academic Publishers, Netherlands.
\bibitem{A2}  K. ADJAMAGBO, Sur les morphismes injectifs et les isomorphismes des vari\'et\'es
alg\'ebriques affines, Communications in Algebra, 24 (3), 1117-1123 (1996).
\bibitem{A3}  K. ADJAMAGBO, On isomorphisms of factorial domains and the Jacobian Conjecture in any Characteristic, Pr\'epublication 91, Octobre 1996, Institut de Math\'ematiques de Jussieu, Universit\'e Paris 6 et Universit\'e Paris 7.
\bibitem{A-E} K. ADJAMAGBO, A. van den ESSEN, Eulerian systems of partial differential equations and the Jacobian Conjecture, Journal of Pure and Applied Algebra 74 (1991) 1-15.
\bibitem{A-K} A. ALTMAN, S. KLEIMAN, Introduction to Grothendieck Duality Theory, Lecture Notes in Mathematics 146, Springer-Verlag, Berlin, 1970.
\bibitem{B}   H. BASS, Differential Structure of \'etale Extensions of Polynomial Algebras, in
Commutative Algebra, Proceedings of a Microprogram Held, June 15-July 2, 1987,M. Hoster, C.
Hunecke, J.D. Sally, ed., Springer-Verlag, New York, 1989.
\bibitem{BCW} H. BASS, E. CONNELL, D. WRIGHT, The Jacobian Conjecture: Reduction of Degree and formal Expansion of the Inverse, Bull. AMS 7(1982), 287-330
\bibitem{D-G} J. DIEUDONNE, A. GROTHENDIECK, Crit\`eres diff\'erentielles de r\'egularit\'e pour les localis\'es des alg\`ebres analytiques, J. Algebra 5, 1967, 305-324.
\bibitem{EGA IV3} A. GROTHENDIECK, J. DIEUDONNE, El\'ements de G\'eom\'etrie Alg\'ebrique, IV, Troisi\`eme Partie, Publications Mathématiques N° 28, 1966, Institut des Hautes Etudes Scientifiques.
\bibitem{EGA IV4} A. GROTHENDIECK, J. DIEUDONNE, El\'ements de G\'eom\'etrie Alg\'ebrique, IV, Quatri\`eme Partie, Publications Mathématiques N° 32, 1967, Institut des Hautes Etudes Scientifiques.
\bibitem{F}  T. FUJITA, On the topology of non-complete surfaces, J. Fac. Sci. Univ. Tokyo 29
(1982) 503-566.
\bibitem{FI}  G. FISCHER, Complex Analytic Geometry, Lecture Notes in Mathematics, Vol. 538, Springer-verlag, Berlin-New York, 1976.
\bibitem{G}  R. V. GURJAR, Affine varieties dominated by ${\Bbb C}^2$, Comment. Math. Helvetici 55(1980) 378-389. 
\bibitem{G-M} R. V. GURJAR, M. MIYANISHI, Affine surface with $\overline{\kappa} \leq \; 1$, in
Algebraic Geometry and Commutative Algebra in Honour of Misayoshi NAGATA, 1987, 99-124.
\bibitem{G-S} Alexander GROTHENDIECK, Hamet SEYDI, Platitude d'une adh\'erence sch\'ematique et lemme de Hironaka g\'en\'eralis\'e, Manuscripta Math. 5 (1971), 323-339.
\bibitem{H}   R. HARTSHORNE, Algebraic Geometry, G.T.M. 52, Springer-Verlag, 1977.
\bibitem{H'}  R. HARTSHORNE, Ample Subvarieties of Algebraic Varieties, Lecture Notes in Math. 156, Springer-verlag, Berlin, 1970.
\bibitem{K}   V. KULIKOV, Generalised and local jacobian problems, Russian Acad. Sci. Izv. Math.
Vol. 41, 1993, No. 2, 351-365.
\bibitem{L}   S. LANG, Algebra, Addison-Wesley, 8th Printing, 1985.
\bibitem{M}  M. MIYANISHI, Non-complete Algebraic Surfaces, Lecture Notes in Math. 857,
Springer-Verlag, Berlin, 1981.
\bibitem{N}  M. NAGATA, Field Theory, Marcel Dekker, Inc., New York,1977.
\bibitem{O}  S. YU. OREVKOV, On three-sheeted polynomial mappings of ${\Bbb C}^2$, Math. URSS Izv., Vol. 29, 1897, No. 3, 587-596.
\bibitem{OD1}  S. ODA, The Last Aproach to the Settlement of the Jacobian Conjecture, arXiv:math/0307080, 2003-2007.
\bibitem{OD2}  S. ODA,  A Valuation Theoretic Approach to the Jacobian Conjecture, arXiv:0706.1138, 2007-2011.
\bibitem{OO}  F. OORT, Units in number fields and function fields, Exposition. Math. 17 (1999), no. 2, 97-115.
\bibitem{R}  C. P. RAMANUJAM, A topological characterisation of the affine plane as an algebraic
variety, Ann. of Math. 94 (1971) 69-88.
\bibitem{SGA} A. GROTHENDIECK, S\'eminaire de G\'eom\'etrie Alg\'ebrique, 1960-1961, I.H.E.S., 1960, Fascicule 1.
\bibitem{S1} H. SEYDI, La Conjecture Jacobienne, Rend. Sem. Mat. Messina, dere II, Vol. VI (1999), 175-179.
\bibitem{S2} H. SEYDI, La Conjecture Jacobienne II, Afr. diaspora J. Math. (N.S.) 10 (2010), no. 1, 87-121.
\bibitem{Z}   M. G. ZAIDENBERG, On Ramanujam surfaces, ${\Bbb C}^{**}$-families, and exotic algebraic structures on ${\Bbb C}^n$, Trans. Moscow Math. Soc., 1994.

\end{thebibliography}
\end{document}